\documentclass[12pt,leqno]{article}

\catcode`\@=11
\long\def\@makecaption#1#2{
    \vskip 10pt
    \setbox\@tempboxa\hbox{\bf #1: \sf #2}
    \ifdim \wd\@tempboxa >\hsize \bf #1: \sf #2\par \else \hbox
to\hsize{\hfil\box\@tempboxa\hfil}
    \fi}
\catcode`\@=12
\usepackage{amsmath}
\usepackage{amsfonts}

\def\kkk{\null\hfill $\Box $ \\}
\def\r{ \right}
\def\eps{\varepsilon}

\def\no{\|\cdot\|}

\newcommand{\km}{k\mbox{-}\min}

\newcommand{\kjm}{(k-j+1)\mbox{-}\min}

\newtheorem{theo}{Theorem}
\newtheorem{lem}[theo]{Lemma}
\newtheorem{prop}[theo]{Proposition}

\newtheorem{cl}[theo]{Claim}

\usepackage{epsfig}
\usepackage{latexsym}
\def\kkk{\null\hfill $\Box $ \\}

\begin{document}
\title{Uniform estimates for order statistics
and Orlicz functions
}

\author{Yehoram Gordon
\thanks{Partially supported by AIM, Palo Alto}
\thanks{Supported in part by ``France-Israel Cooperation Grant \#3-1350"
and by the ``Fund for the Promotion of Research at the Technion".
}
\and
Alexander Litvak $\mbox{}^{\dagger}$
\and
Carsten Sch\"utt$\mbox{}^{\dagger}$
\and
Elisabeth Werner$\mbox{}^{\dagger}$ \thanks{Partially supported by a NSF
Grant, by a FRG Grant, and by a BSF Grant}
}
\date{}
\maketitle

\begin{abstract}
We establish uniform estimates for order statistics: Given a 
sequence of independent identically distributed random variables 
$\xi_{1},\dots,\xi_{n}$ with log-concave distribution and scalars 
$x_{1},\dots,x_{n}$, for every $k\le n$ we provide estimates for 
$\mathbb E \, \, \km_{1\leq i\leq n} |x_{i}\xi _{i}|$ and 
$\mathbb E \operatornamewithlimits{k- max}_{1\leq i\leq n}|x_{i}\xi_{i}|$ 
in terms of the value $k$ and the appropriate Orlicz norm 
$\|(1/x_{1},\dots,1/x_{n})\|_M$, associated with the distribution 
function of the random variable $|\xi_{1}|$. For example, if $\xi_1$ is the 
standard $N(0,1)$ Gaussian random variable, then the corresponding 
Orlicz function is $M(s)=\sqrt{\tfrac{2}{\pi}}\int_{0}^{s} e^{-\frac{1}{2t^{2}}}dt$. 
We would like to emphasize that our estimates do not depend on the length $n$ of the 
sequence. 
\end{abstract}

\section{Introduction}
In this paper we establish uniform estimates for order statistics.
The $k$-th order statistic of a statistical sample  of size $n$ is
equal to its $k$-th smallest value, or equivalently its $(n-k+1)$-th
largest value. Order statistics are among the most fundamental tools
in non-parametric statistics and inference and consequently there is
extensive literature on order statistics. We only cite \cite{AN, BCo, 
BCS, BR1, DN, W1} and references therein.

Order statistics are more resilient to faulty sensor reading than
max, min or average and thus they find applications when  methods
are needed to study configurations that take on a ranked order. To
name only a few: Wireless networks, signal processing, image
processing, compressed sensing, data reconstruction, learning theory
and data mining. A sample of works done in this area are
\cite{BCh, BaDaDeVWa, BL, BR2, CRT, CoDDeV, DG, D, GLMP, MZ, MPT}.

Order statistics on random sequences appear naturally in Banach space theory, 
in computations of  various random parameters
associated with the geometry of convex bodies in high dimensions, 
in random matrix theory (computing the distribution of eigenvalues), 
and in approximation theory (see e.g. \cite{Go2, GGMP, GLMP, GLSW2, 
Gu, MaPa, RU, SZ1, SZ2}). 
This list of course does not include the enormous quantity of
published works which deal with evaluations and applications of
$max$ and $min$ associated with various random parameters, e.g.,
smallest and largest eigenvalues of random matrices, as these are
the extreme values in the scale of order statistics.

For these important special cases of  order statistics,  the minimum and
maximum value of a sample, very precise estimates were obtained in 
\cite{GLSW, GLSW3, GLSW4}.  The new approach
started  there was to give  estimates of the  minimum and maximum value of the sample
\begin{equation}\label{minmax}
\mathbb E\min_{1\leq i\leq n}|x_{i} \xi_{i}| \ \ \  \mbox {and} \ \ \
\mathbb E\max_{1\leq i\leq n}|x_{i} \xi_{i}|,
\end{equation}
in terms of {\em Orlicz norms} (see the definition below).
The expressions for the estimate in  case (\ref{minmax}) are
relatively simple. For instance, it was shown in \cite{GLSW} that
$$
c_{1}\|x\|_{M}  \leq \mathbb{E}
\max_{1\leq i\leq n}|x_{i}\xi_{i}(\omega)|
\leq c_{2}\|x\|_{M},
$$
and in \cite{GLSW3, GLSW4} that 
$$
c_3 \,
\left(\sum_{i=1}^{n}\frac{1}{|x_{i}|}\right)^{-1}
\leq\mathbb E\min_{1\leq i\leq n}|x_{i} \xi_{i}|
\leq c_4\,
\left(\sum_{i=1}^{n}\frac{1}{|x_{i}|}\right)^{-1},
$$
where $c_1$, $c_2$, $c_3$, and $c_4$  are absolute positive constants
and $\|\cdot \|_{M}$ is an Orlicz norm, depending on the distribution of 
$\xi_1$ only.
In fact, in \cite{GLSW} much more general case was considered 
(see also \cite{J} and \cite{MS}). 

\vskip 3mm
Here we study the values 
\begin{equation}\label{general}
\mathbb E \, \, \km_{1\leq i\leq n} |x_{i}\xi _{i}|  \ \ \  \mbox {and} \ \ \
\mathbb E
  \operatornamewithlimits{k- max}_{1\leq i\leq n}|x_{i}\xi_{i}|,
\end{equation}
for general order statistics for i.i.d. (independent identically distributed) 
random variables $\xi_{1},\dots,\xi_{n}$ and scalars $x_{1},\dots,x_{n}$, 
where for a given sequence of real numbers $a_{1},\dots,a_{n}$ we denote
the $k$-th smallest one by $\km _{1\leq i\leq n}a_{i}$. In particular,
$1\mbox{-}\min _{1\leq i\leq n}a_{i}  = \min _{1\leq i\leq n}a_{i}$ and
$n\mbox{-}\min _{1\leq i\leq n}a_{i} =  \max _{1\leq i\leq n}a_{i}$.
In the same way we denote the $k$-th biggest number by
$k\mbox{-}\max_{1\leq i\leq n}a_{i}$. Thus,
$k\mbox{-}\max_{1\leq i\leq n}a_{i} =(n-k+1)\mbox{-}\min _{1\leq i\leq n}a_{i}$.
In fact, in the theory of order statistics the standard notation for $\km$ is  
$a_{k:n}$. In this paper such a notation could be misleading and we prefer 
to use $\km$. 
\par
Now the expressions  get more involved than in  case (\ref{minmax}).
In view of possible applications we strive to keep them
as simple as  possible -- at the expense of the constants involved.
We show that if $\xi_1$ has a log-concave distribution then for $1\leq k\leq n/2$
\begin{eqnarray*}
c_{1} \ \max_{1 \leq j \leq k}\   \| \left(1/x_i\r)
_{i=j}^n \|_{\frac{2e}{k-j+1}N}^{-1}
\leq \mathbb E \, \, \km_{1\leq i\leq n} |x_{i}\xi _{i}|
\leq  c_{2} \max_{1 \leq j \leq k}\
\|\left(1/x_i\r)_{i=j}^n\|^{-1}_{\frac{2e}{k-j+1}N},
\end{eqnarray*} 
and for $1\leq k\leq cn$
\begin{eqnarray*}
&&
c_{3} \left( \max_{0 \leq \ell \leq ck-1}\   \| \left(1/x_i
  \r)_{i=1}^{k+\ell} \|_{\frac{2e}{\ell+1}N}^{-1} +
  \|(x_{k+ck},\dots,x_{n})\|_{M}\right)
  \\&&
  \leq\mathbb E
  \operatornamewithlimits{k-max}_{1\leq i\leq n}|x_{i}\xi_{i}|
 \\&&
\leq c_{4}\left( \max_{0 \leq \ell \leq ck-1}\
  \|\left(1/x_i\r)_{i=1}^{k+\ell}\|^{-1}_{\frac{2e}{\ell+1}N}
+  \|(x_{k+ck},\dots,x_{n})\|_{M}\right),
\end{eqnarray*}
where $N$, $M$ are {\it Orlicz functions}  (see the definitions below) 
and $\|\cdot \|_{N}$, $\|\cdot \|_{M}$ are the corresponding Orlicz norms.
The Orlicz functions $N$, $M$ are computed in terms of the distribution function
of the random variables under consideration. The constants $c$, $c_1$, $c_{2}$, 
$c_3$, $c_{4}$ depend -- mildly -- on the distribution function of the random variables 
and on $k$ (of the order of $\ln k$ or $1/\ln k$), but -- and this  is the important point --
they do not depend in any way on the number $n$ and on the scalars
$x_{1},\dots,x_{n}$. The precise statements are given in Section~\ref{mainres}.
We would like to note that Orlicz functions appear naturally in the connection with 
log-concave distributions. For example in the important work of Gluskin and Kwapie\'n 
\cite{GK} Orlicz functions were used to obtain tail and moment estimates for 
sums of independent random variables. Recently, Latala \cite{Lat} proved tail 
comparison theorem for log-concave vectors.

In problems where only a small number of random variables is involved, 
numerical computations will give sufficient estimates for order statistics.
However, in the case when a large number of random variables is involved,
numerical computations may not be feasible. Our formulae allow easy computations 
also in that situation.

Finally let us mention that throughout this paper we use the following 
notation. For a random variable $\xi $ on a probability space
$(\Omega,\mathfrak A,\mathbb P)$ we denote its distribution function
by $G_{\xi}$ and $1-G_{\xi}$ by $F_{\xi}$
$$
G_{\xi}(t) =\mathbb P \left(\left\{ \xi \leq t\r\}\r)
\hskip 15mm\mbox{and}\hskip 15mm
F_{\xi} (t) =\mathbb P \left(\left\{ \xi > t\r\}\right) .
$$

\medskip
\noindent
{\bf Acknowledgment.} 
We would like to thank Hermann K\"onig, Kiel, for discussions.

\section{Preliminaries. Orlicz functions and norms.}
\label{secorl}

In this section we recall some facts about Orlicz functions and norms. 
For more details and other  properties of Orlicz spaces we refer to \cite{KrRu, LT, RaRe}.
\par

A left continuous convex function $M:[0,\infty)\rightarrow [0,\infty] $ is called
Orlicz function or Young function, if $M(0)=0$ and if $M$ is neither
the function that is constant $0$ nor the function that takes the value $0$ at $0$
and is $\infty$ elsewhere.
The corresponding Orlicz norm on $\mathbb R^{n}$ is defined by
\begin{equation}\label{Functional}
 \left\| x \right\| _M=
\inf \left\{ \rho > 0\left| \, \,  \sum_{i=1}^n
 M\left( |x_i |/\rho\right) \leq 1 \right.\right\} .
\end{equation}
\par
Note that the expression for $\|\cdot\|_M$ makes also sense if the function $M$ is
merely positive and increasing. Although in that case the expression need not be a
norm, we keep the same notation $\|\cdot\| _M$. We often use  formula (\ref{Functional}) in a
slightly different form, namely
$$
  1/ \left\| x \right\| _M  = \sup \left\{ \rho > 0\left| \, \,
  \sum_{i=1}^n M\left(\rho\ |x_i |\right) \leq 1 \right.\right\}.
$$
Clearly, $M\leq \bar M$ implies $\no _M \leq \no _{\bar M}$.  Moreover, if $M$ is
an Orlicz function and $s\geq 1$,  then  
\begin{equation} \label{conin}
         s M(t) \leq M(st)  
\end{equation} 
for every $t\geq 0$. In particular, this  implies 
\begin{equation} \label{orlcomp}
\|\cdot \| _{sM} \leq s \|\cdot \| _M.
\end{equation}

The dual function $M^{*}$ to an Orlicz function $M$ is defined
by
$$
  M^{*}(s)=\sup_{0\leq t<\infty}(t\cdot s-M(t)).
$$
For instance,  for  $M(t)=\frac{1}{q}\ t^{q}$, $ q \geq 1$,  the dual function is
$M^{*}(t)=\frac{1}{q^{*}}\ t^{q^{*}}$ with $\frac{1}{q}+\frac{1}{q^{*}}=1$.
\par
Let the function $p=p_M : [0,\infty)\rightarrow[0,\infty]$ be given by
$$
 p(t)=\left\{
 \begin{array}{cc}
 0   &  t=0   \\
 M^{\prime}(t)   &   M(t)<\infty   \\
 \infty   &   M(t)=\infty,
 \end{array}
 \right.
$$
where $M^{\prime}$ is the left hand side derivative of $M$. Then  $p$ is increasing and the
left hand side inverse $q$ of the increasing function $p$ is
$$
  q(s)=\inf\{t\in[0,\infty)\ |\ p(t)>s\}.
$$
Then
$$
  M^{*}(s)=\int_{0}^{s}q(t)dt.
$$
To a given random variable $\xi$ we associate an Orlicz function
$M=M_{\xi}$ in the following way: 
\begin{equation}\label{Orlicz1}
 M(s)=\int_{0}^{s}\int_{\frac{1}{t}\leq|\xi|}|\xi|d\mathbb P dt
 = \int_{\frac{1}{s}\leq|\xi|} \left(s |\xi| - 1\right)  d\mathbb P.
\end{equation}
The equality here follows by changing the order of integration and the 
convexity of $M$ follows by the definition of convexity. We prefer to keep 
in mind both formulae for $M$. Note that equivalently one can write 
$$
   M(s) = \mathbb E \left(s|\xi| - 1 \right)_{+}, 
$$
where, as usual, $h_+(x)$ denotes $h(x)$ if $h(x)\geq 0$ and $0$ otherwise.

\medskip

We claim that the dual function $M^{*}=M^*_{\xi}$ is given on
$[0,\int_{}^{}|\xi_{}|d\mathbb P]$ by
\begin{equation}\label{OrliczMin}
M^{*}\left(\int_{t\leq|\xi|}|\xi_{}|d\mathbb P\right)
=\mathbb P(|\xi|\geq t)
\end{equation}
and $M^{*}(s)=\infty$ for $s>\int_{}^{}|\xi_{}|d\mathbb P$.

Indeed, by  definition 
$$
  M^{*}(s)=\sup_{0\leq w}(w\cdot s-M(w)) =\sup_{0\leq w} \left( w\cdot s
  -\int_{0}^{w}\int_{\frac{1}{u}\leq|\xi|}|\xi|d\mathbb P du\right) 
$$
$$
  = \sup_{0\leq w} \int_{0}^{w}\left(s- \int_{\frac{1}{u}\leq|\xi|}|\xi|d\mathbb 
  P\right) du. 
$$
If $s>\int_{}^{}|\xi_{}|d\mathbb P$ then the supremum is equal to $\infty$. 
Now fix $t\geq 0$, set 
$$
   s = \int_{t\leq|\xi|}|\xi_{}|d\mathbb P 
$$
and consider the function
$$
  \phi(w):= \int_{0}^{w}\left(s - \int_{\frac{1}{u}\leq|\xi|}|\xi|d\mathbb 
  P\right) du .
$$
It is easy to see that $\phi$ is increasing on $[0, 1/t]$ and decreasing 
on $[1/t, \infty)$. Therefore,
$$
  M^*(s) = \sup_{0\leq w}  \phi(w) =  \phi(1/t) = \int_{0}^{1/t} \, 
  \int_{t\leq |\xi| < 1/u} \ |\xi| \ d\mathbb P \ du . 
$$
Changing the order of integration we obtain 
$$
  M^*(s) = \int_{t\leq |\xi| } \ \int_{0}^{1/|\xi| } \, \ |\xi| \ du \ 
  d\mathbb P =  \int_{t\leq |\xi|} \ d\mathbb P = \mathbb P(|\xi|\geq t), 
$$
which proves (\ref{OrliczMin}).

\bigskip

\bigskip

In the Gaussian case we have
$$
  F(t)=\mathbb P \left(\left\{| \xi |> t\r\}\right)=
  \sqrt{\tfrac{2}{\pi}} \int_{t}^{\infty}e^{-\frac{s^{2}}{2}}ds
$$
and thus
\begin{equation}\label{OrlGauss}
M(s)
=\sqrt{\tfrac{2}{\pi}}\int_{0}^{s}
\int_{\frac{1}{t}}^{\infty}ue^{-\frac{u^{2}}{2}}dudt
=\sqrt{\tfrac{2}{\pi}}\int_{0}^{s}e^{-\frac{1}{2t^{2}}}dt .
\end{equation}
This implies that  on the interval $[0,\sqrt{2/\pi}]$ $M^{*}$ is given by
$$
M^{*}(s)=\int_{0}^{s}
\frac{1}{\sqrt{2\ln\left(\sqrt{\tfrac{2}{\pi}}\frac{1}{u}\right)}}du.
$$
For $s>\sqrt{2/\pi}$,   $M^{*}(s)=\infty$.
\vskip 3mm

\section{The main results}\label{mainres}

Now we consider certain functions associated with
a random variable $\xi:\Omega\rightarrow \mathbb R$. 
\par
The function
$F:[0,\infty)\rightarrow[0,\infty)$
is given by
\begin{equation}\label{distrrest}
F(t)
=\mathbb P(|\xi|>t) .
\end{equation}
We assume that $F$ is strictly decreasing on $[0,\infty)$ and
$F(0)=1$. In particular, $F$ is invertible.
\par
The function $N:[0,\infty)\rightarrow[0,\infty)$ is defined by
\begin{equation}\label{OrliczDistr}
N(t)=\ln\frac{1}{F(t)}  
\end{equation}
and is assumed to be convex. In particular, $N$ is
an Orlicz function. For such a function $N$ and $k\in\mathbb N$ we put
\begin{equation}\label{OrlNj}
N_{j}=\frac{2e}{k-j+1}N,
\hskip 20mm
j=1,\dots,k .
\end{equation}
\par
Furthermore, let us observe that under assumptions above for all $t\geq0$
and all $s\geq1$ we have
\begin{equation}\label{ExpFunc}
F(s t)\leq F(t)^{s} .
\end{equation}
Indeed, by (\ref{conin}) we have $s N(t)\leq N(st)$, i.e.
$- s \ln F(t) \leq - \ln F(st)$, 
which is equivalent to (\ref{ExpFunc}).

The following theorem generalizes results from \cite{GLSW3, GLSW4}, where
similar estimates were obtained for Gaussian distributions. Of course, the
Gaussian case is simpler and the corresponding formulae are less involved.
We discuss the details in Remark 1 after the theorem.

\begin{theo}\label{Kmin1}
Let $1\leq k\leq\frac{n}{2}$ and
let $\xi_{1},\dots,\xi_{n}$ be i.i.d. copies of a random variable $\xi$.
Let $F$, $N$ and $N_{j}$, $j=1,\dots,k$, be
as specified in (\ref{distrrest}), (\ref{OrliczDistr}) and (\ref{OrlNj}).
Then for all $0< x_1\leq  x_2\leq \dots \leq x_n$ 
\begin{eqnarray*}
c_{1}
\ \max_{1 \leq j \leq k}\   \| \left(1/x_i\r) _{i=j}^n \|_{N_{j}}^{-1}
\leq \mathbb E \, \, \km_{1\leq i\leq n} |x_{i}\xi _{i}|
\leq 16 e^2 \ C_{N}\ \ln(k+1)
\max_{1 \leq j \leq k}\   \|\left(1/x_i\r)_{i=j}^n\|^{-1}_{N_{j}},
\end{eqnarray*}
where $c_{1}=1- \tfrac{1}{\sqrt{2 \pi}}$ and $C_{N}=\max\{N(1),1/N(1)\}$.

Moreover, the lower estimate does not require the condition
``$N$ is an Orlicz function".
\end{theo}
\vskip 3mm

\noindent
{\bf Remark (the Gaussian case).}
In \cite{GLSW3, GLSW4} it was shown that for $N(0,1)$ random variables
$g_{i}$, $i=1,\dots,n$
and for all
$0< x_1 \leq x_2 \leq \dots \leq x_n$
\begin{eqnarray}\label{old}
c_0
\max_{1 \leq j \leq k}
\ \frac {k+1-j}{\sum_{i=j}^n 1/x_i}
\leq \mathbb E\, \,
\km_{1\leq i\leq n} |x_{i}g_{i}| \leq
2\sqrt{2 \pi}\, \ln(k+1)\, \max_{1 \leq j \leq k}
\ \frac {k+1-j}{\sum_{i=j}^n 1/x_i }.
\end{eqnarray}
where
$c_0
=\left(1-\frac{1}{\sqrt{2\pi}}\right)\frac{1}{2 e}
\sqrt{\frac{\pi}{2}}$.
In Section~\ref{gausscl} we show that the Gaussian distribution satisfies
the conditions of Theorem~\ref{Kmin1}. Thus the estimate (\ref{old}) can be obtained
from Theorem~\ref{Kmin1} (with different absolute constants). 

\bigskip

Our second theorem provides bounds for expectations of $k$-max. 
As in Theorem \ref{Kmin1} we assume that $F$ is strictly decreasing, $F(0)=1$,  
and that $N=-\ln F$ is a convex function, where $F$ is given by (\ref{distrrest}). 
Note that such a function $F$ satisfies
\begin{equation}\label{SubMultEst3}
\int_{t\leq|\xi_{1}|}|\xi_{1}|d\mathbb P\leq \left(1+ \frac{1}{N(t)} \r) t\cdot F(t) 
\end{equation} 
for all positive $t$. 
We verify this. Since $F=e^{-N}$ and $N$ is convex 
$$
\int_{t\leq|\xi_{1}|}|\xi_{1}|d\mathbb P
=t\cdot F(t)+\int_{t}^{\infty}F(s)ds
=t\cdot F(t)+\int_{t}^{\infty}e^{-N(s)}ds . 
$$
Using (\ref{conin}), we have $N(s)\geq \frac{s}{t}N(t)$ for $s\geq t$. 
Therefore
\begin{eqnarray*}
\int_{t\leq|\xi_{1}|}|\xi_{1}|d\mathbb P
&\leq& t\cdot F(t)+\int_{t}^{\infty}e^{-\frac{s}{t}N(t)}ds   \\
&\leq& t\cdot F(t)+\frac{t}{N(t)}e^{-N(t)}
= t\cdot F(t)+\frac{t}{N(t)}F(t) , 
\end{eqnarray*}
which  implies (\ref{SubMultEst3}).

\medskip

\begin{theo}\label{TheoK-max}
Let $\xi_{1},\dots,\xi_{n}$ be i.i.d. copies of a random variable $\xi$.
Let $F$, $M$, and $N$ be as specified in (\ref{distrrest}),
(\ref{Orlicz1}), and (\ref{OrliczDistr}).
Let $1 < k \leq n$ and $k_{0}=\left[ \frac{4(k-1)}{F(1)}\right]$.
Assume that $k+k_0\leq n$. Then for all
$x_{1}\geq x_{2}\geq\cdots\geq x_{n}>0$ 
$$ 
 \frac{1}{4}\ \left( 
 \max_{0 \leq \ell \leq k_{0}-1}\   \| \left(1/x_i
  \r)_{i=1}^{k+\ell} \|_{\frac{2e}{\ell+1}N}^{-1} 
 + \left( 1 + \frac{\ln(8(k-1))}{N(1)}\right)^{-1} \
  \|(x_{k+k_{0}},\dots,x_{n})\|_{M} \r) 
$$
$$
 \leq\mathbb E
 \operatornamewithlimits{k-max}_{1\leq i\leq n}|x_{i}\xi_{i}|
 \\
 \leq c \left(C_{N}\ln(k+1) \max_{0 \leq \ell \leq k_{0}-1}\
 \|\left(1/x_i\r)_{i=1}^{k+\ell}\|^{-1}_{\frac{2e}{\ell+1}N}
 + \ \|(x_{k+k_{0}},\dots,x_{n})\|_{M} \r),
$$ 
where $C_N = \max\{N(1),1/N(1)\}$, and  
$c$ is an absolute positive constant.
\end{theo}
\vskip 3mm

\noindent
{\bf Remark.} The case $k=1$ was obtained in \cite{GLSW} (see also 
Lemma~\ref{1max} below): Let
$$
  M(s) =\int_{0}^{s}\int_{\frac{1}{t}\leq|\xi_1|}|\xi_1|d\mathbb P dt
  = \int_{\frac{1}{s}\leq|\xi|} \left(s |\xi| - 1\right)  d\mathbb P.
$$
Then, for all $x\in\mathbb R^{n}$ one has
$$
  c_{1}\|x\|_{M} \leq\int_{\Omega}\max_{1\leq i\leq n}|x_{i}\xi_{i}(\omega)|
  d\Bbb P(\omega) \leq c_{2}\|x\|_{M},
$$
In particular, in the Gaussian case (\ref{OrlGauss}),
$$
M(s)
=\sqrt{\tfrac{2}{\pi}}\int_{0}^{s}e^{-\frac{1}{2t^{2}}}dt .
$$

\section{k-min}

We need the  following two simple lemmas. Similar lemmas were used in
\cite{GLSW3, GLSW4}. For the sake of completeness we  provide the proofs.

\begin{lem}\label{lemmin}
Let $0<x_1\leq x_2\leq ... \leq x_n$.
Let  $\xi _1, \dots, \xi_n$ be i. d. random variables.
Let $F(t)=\mathbb P\{|\xi_{1}|> t\}$ and $G(t)=1-F(t)$.
Then
$$
\mathbb P  \left\{ \min_{1\leq i\leq n} |x_{i}\xi_{i}|
\leq t \right\}
\leq \sum _{i=1}^{n} G\left(t/ x_i\r).
$$
Moreover, if the $\xi _{i}$'s are independent then for every
$t>0$
$$
\mathbb P\left\{\min_{1\leq i\leq n}|x_{i} \xi _{i}|> t \right\}
=\prod _{i=1}^{n} F\left(t/ x_i\r).
$$
\end{lem}

\medskip

\noindent{\bf Proof. } 
 Denote $A_{k}(t)=\{\omega \ |\ |x_{k} \xi_{k}(\omega)|> t\} =
\{\omega \ |\ | \xi_{k}(\omega)|> t/x_k \} $ and
$$
  A(t) = \{\omega \ | \ \min_{k\leq n}|x_{k} \xi_{k}(\omega)| >t\}
  = \bigcap _{k\leq n} A_k (t).
$$
 Then
$$
  \mathbb P \left( A(t) \r) \geq 1-\sum_{k=1}^{n} \mathbb
  P\left( A_{k}(t)^{c} \r) = 1- \sum_{k=1}^{n}\ G\left(t/x_{k}\r) ,
$$
which proves the first estimate.
The second estimate is trivial.
\kkk

For the second lemma we need the following Proposition,
proved in \cite{GLSW3}.

\begin{prop}\label{agmean}
Let $1\leq k\leq n$. Let $a_{i}$, $i=1,\dots,n$ be nonnegative
real numbers. Assume
$$
   0<  a:= \frac{e}{k}\, \sum _{i=1}^n a_i < 1 .
$$
Then
$$
 \sum _{l=k}^n
\sum_{A\subset \{1, 2, ..., n\} \atop |A|=l}\prod_{i\in A}a_{i}
<
\frac{1}{\sqrt{2 \pi k}} \, \,  \frac{a^k}{1-a}.
$$
\end{prop}
\vskip 3mm

\begin{lem}\label{forthtwo}
Let $1\leq k\leq n$.
Let $0<x_1\leq x_2\leq ... \leq x_n$ and $\xi _1, \dots, \xi_n$ be
i.i.d. random variables.
Let $G(t) =\mathbb P\{|\xi_{1}|\leq t\}$ and
$$
   a = a(t) = \frac{e}{k}\ \sum _{i=1}^{n} G\left(t/ x_i\r).
$$
 Assume that $t$ is such that  $0<a<1$. Then
\begin{equation} \label{newkminn}
  \mathbb P \left\{\km_{1\leq i\leq n} |x_{i} \xi _{i}|\leq
t\right\}
 \leq   \frac{1}{\sqrt{2 \pi k}} \, \, \frac{a^k}{1-a}.
\end{equation}
\end{lem}

\medskip

\noindent
{\bf Remark.} Note that if $G$ is continuous and $G(s)= 0$
if and only if $s=0$ then the condition on $t$ in Lemma \ref{forthtwo}
above corresponds to the condition
$$
  0<t< \| (1/x_i) _{i=1}^{n} \|_H^{-1},
$$
where $H=\frac{e}{k}G$.

\medskip
\noindent
{\bf Proof of Lemma~\ref{forthtwo}.  }
We have
\begin{eqnarray*}
&&\mathbb P\left\{\omega\left|\
\km_{1\leq i\leq n}|x_{i} \xi_{i}(\omega)|\leq
 t \right.\right\}  \\
&&= \mathbb P\left\{\omega\left|\exists i_{1},\dots,i_{k}
\geq 1 : | \xi_{i_{j}}(\omega)|\leq\frac{ t }{x_{i_{j}}}
\right.\right\}   \\
&&= \mathbb P\left(\bigcup_{\ell=k}^{n}
\bigcup_{A\subset \{1,...,n\}
\atop |A|=\ell} \left\{\omega\left|\forall i\in A:|
\xi_{i}(\omega)| \leq
\frac{ t }{x_{i}}\hskip 1mm\mbox{and}\hskip 1mm
\forall i\notin A:| \xi_{i}(\omega)|>\frac{ t
}{x_{i}}\right.\right\} \right)  \\
&&=
\sum_{l=k}^{n} \ \sum_{A\subset \{ 1, ..., n\} \atop |A|=l}\
 \prod_{i \in A}
\mathbb P \left\{\omega \ \left| \  | \xi_{i}(\omega)| \leq
\frac{ t }{x_{i}} \right.\right\} \,
\prod_{i \notin A}
\mathbb P \left\{\omega \left| \ |\xi_{i}(\omega)|>\frac{ t
}{x_{i}}\right.\right\} .
\end{eqnarray*}
It follows
\begin{eqnarray*}
\mathbb P\left\{\omega\left|\
\km_{1\leq i\leq n}|x_{i} \xi_{i}(\omega)|\leq
 t \right.\right\}
 &\leq &
\sum_{l=k}^{n} \ \sum_{A\subset \{ 1, ..., n\} \atop |A|=l}\
 \prod_{i \in A} \mathbb P \left\{\omega \ \left| \  |
\xi_{i}(\omega)| \leq
\frac{ t }{x_{i}} \right.\right\}    \\
&\leq&
\sum_{l=k}^{n} \ \sum_{A\subset \{ 1, ..., n\} \atop |A|=l}\
\prod_{i \in A}  G\left( t/x_i\r) .
\end{eqnarray*}
Proposition~\ref{agmean} implies the desired result.
\kkk

\begin{lem}\label{Partition}
Let
$H:\mathbb R\rightarrow\mathbb R$ be an Orlicz function.
For every $k$ with $1\leq k\leq n$ and every
$0< x_{1}\leq x_{2}\leq\cdots\leq x_{n}$ there is a partition of nonempty sets
$A_{1},\dots,A_{k}$ of the set $\{1,\dots,n\}$ such that
\begin{equation}\label{Partition3}
\min_{1\leq j\leq k}
\left\|\left(\frac{1}{x_{i}}\right)_{i=j}^{n}\right\|_{\frac{H}{k-j+1}}
\leq
4 \max\{H(1),1/H(1)\}\min_{1\leq j\leq k}
\left\|\left(\frac{1}{x_{i}}\right)_{i\in A_{j}}\right\|_{H} .
\end{equation}
\end{lem}
\vskip 2mm

We want to emphasize that it is important that the partition consists of exactly
$k$ sets. Our proof shows that the partition
can be taken as intervals, that is $A_j=\{n_{j}+1, \ldots,  n_{j+1}\}$ for an
increasing sequence $0=n_0 < n_1 < \ldots < n_{k}=n$.
\vskip 2mm

\noindent
{\bf Proof.} We may assume that $H(1)=1$.
Indeed, as $H$ is convex and as $H(0)=0$, 
$H(s)\leq\frac{s}{t}H(t)$ for all $0<s<t$. Thus
if $H(1)\leq 1$, then
$$
  H(1)\|y\|_{\frac{H}{H(1)}} \leq\|y\|_{H}\leq\|y\|_{\frac{H}{H(1)}}
$$
for every $y\in\mathbb R^{n}$.
Similarly, if $H(1)>1$
$$
  \|y\|_{\frac{H}{H(1)}}
  \leq\|y\|_{H}\leq H(1)\|y\|_{\frac{H}{H(1)}}.
$$

We consider three cases.

\smallskip

\noindent
{\it Case 1: \ }
\begin{equation}\label{Partition2}
\frac{1}{x_{1}} \leq\frac{1}{4}\left\|\left(\frac{1}{x_{i}}
\right)_{i=1}^{n}\right\|_{\frac{H}{k}} .
\end{equation}

Note that  $H(1)=1$ implies $t=\| (t, 0, \ldots, 0) \|_H$ for every $t>0$,
in particular $1/x_1 = \| (1/x_1, 0, \ldots, 0) \|_H$.
We put $n_{0}=0$ and after having chosen $n_{0},\dots,n_{\ell} <n$ we define
$n_{\ell+1}\leq n$ to be the largest integer such that
\begin{equation}\label{Partition1}
\left\|\left(\frac{1}{x_{i}}\right)_{i=n_{\ell}+1}^{n_{\ell+1}}\right\|_{H}
\leq\frac{1}{2}\left\|\left(\frac{1}{x_{i}}\right)_{i=1}^{n}\right\|_{\frac{H}{k}} .
\end{equation}
We define
$$
B_{\ell}=\{n_{\ell-1}+1,\dots,n_{\ell}\},
\hskip 20mm
\ell=1,\dots,L .
$$
These sets are basically the partition
we are looking for, except for a slight change that is necessary in order to get
exactly $k$ sets.
\par
We verify first that such a partition exists. For this we have
to show that each $B_{\ell}$ contains at least one element,
i.e. $B_{\ell}\ne\emptyset$. In other words, we show that
$0=n_0<n_1<\ldots <n_L=n$. Indeed, if $n_{l-1} <n$, then
$n_{\ell-1}+1\in B_{\ell}$ because
$$
\frac{1}{4}\left\|\left(\frac{1}{x_{i}}
\right)_{i=1}^{n}\right\|_{\frac{H}{k}}
\geq \frac{1}{x_{1}}
\geq \frac{1}{x_{n_{\ell-1}+1}}
=\left\|\left(0,\dots,0,
\frac{1}{x_{n_{\ell-1}+1}},0,\dots,0\right)\right\|_{H} .
$$
In the last equality we used again that $H(1)=1$.
Thus $B_{\ell}\ne\emptyset$ and $n_{L}=n$ which means that the partition
is well defined.
\par
We show now that $L>k$. By (\ref{Partition1}) for every $\eps \in (0, 1)$
and for $\ell=0,\dots,L-1$ we have
$$
  \sum_{i=n_{\ell}+1}^{n_{\ell+1}} H\left((2-\eps) \left\|\left(\frac{1}{x_{i}}
  \right)_{i=1}^{n}\right\|_{\frac{H}{k}}^{-1}\frac{1}{x_{i}}\right) \leq 1,
$$
which implies
$$
  \sum_{i=1}^{n}H\left((2-\eps)\left\|\left(\frac{1}{x_{i}}\right)_{i=1}^{n}
  \right\|_{\frac{H}{k}}^{-1}\frac{1}{x_{i}}\right)\leq L.
$$
Therefore
$$
  \left\|\left(\frac{1}{x_{i}}\right)_{i=1}^{n}\right\|_{\frac{H}{L}}\leq
  \frac{1}{2}\left\|\left(\frac{1}{x_{i}}\right)_{i=1}^{n}\right\|_{\frac{H}{k}}.
$$
This implies $L>k$ and below we use that the  inequality is strict.
\par
We claim that for all $\ell=1,\dots,k$ one has
\begin{equation}\label{Partition3a}
\left\|\left(\frac{1}{x_{i}}\right)_{i\in B_{\ell}}\right\|_{H}
\geq\frac{1}{4}
\left\|\left(\frac{1}{x_{i}}\right)_{i=1}^{n}\right\|_{\frac{H}{k}}.
\end{equation}
Suppose that there is $\ell$ with $1\leq \ell\leq k$ such that
\begin{equation}\label{Partition3b}
\left\|\left(\frac{1}{x_{i}}\right)_{i\in B_{\ell}}\right\|_{H}
<\frac{1}{4}
\left\|\left(\frac{1}{x_{i}}\right)_{i=1}^{n}\right\|_{\frac{H}{k}}.
\end{equation}
Since $L>k\geq \ell$ we have $n_{\ell}+1\leq n$.
$\|\cdot \|_{H}$ is a norm. Therefore, by the triangle
inequality and since $H(1)=1$,
$$
\left\|\left(\frac{1}{x_{i}}
\right)_{i=n_{\ell-1}+1}^{n_{\ell}+1}\right\|_{H}
\leq\left\|\left(\frac{1}{x_{i}}\right)_{i\in B_{\ell}}\right\|_{H}
+\left\|\left(0,\dots,0,\frac{1}{x_{n_{\ell}+1}}\right)\right\|_{H}
=\left\|\left(\frac{1}{x_{i}}\right)_{i\in B_{\ell}}\right\|_{H}
+ \frac{1}{x_{n_{\ell}+1}}.
$$
By (\ref{Partition2}) and (\ref{Partition3b})
$$
\left\|\left(\frac{1}{x_{i}}
\right)_{i=n_{\ell-1}+1}^{n_{\ell}+1}\right\|_{H}
<\frac{1}{2}\left\|\left(\frac{1}{x_{i}}\right)_{i=1}^{n}
\right\|_{\frac{H}{k}} .
$$
This contradicts the definition of $n_{\ell}$.

Now we define the partition $A_{1},\dots,A_{k}$.
We put $A_{\ell}=B_{\ell}$ for $1\leq \ell\leq k-1$ and
$$
  A_{k}=\bigcup_{\ell=k}^{L}B_{\ell}.
$$
Then, by (\ref{Partition3a}),
\begin{eqnarray}\label{partition31}
\min_{1\leq j\leq k}
\left\|\left(\frac{1}{x_{i}}\right)_{i=j}^{n}
\right\|_{\frac{H}{k-j+1}}
&\leq&\left\|\left(\frac{1}{x_{i}}\right)_{i=1}^{n}
\right\|_{\frac{H}{k}} \nonumber  \\
&\leq& 4\min_{1\leq\ell\leq k}
\left\|\left(\frac{1}{x_{i}}\right)_{i\in
B_{\ell}}\right\|_{H}
\leq 4\min_{1\leq\ell\leq k}
\left\|\left(\frac{1}{x_{i}}\right)_{i\in
A_{\ell}}\right\|_{H},
\end{eqnarray}
which proves (\ref{Partition3}).

\smallskip

\noindent
{\it Case 2: \ } 
$$
  \frac{1}{x_{1}} > \frac{1}{4}\left\|\left(\frac{1}{x_{i}}
  \right)_{i=1}^{n}\right\|_{\frac{H}{k}}
\quad \quad \mbox{ and for all } \, \,   j\leq k \, \,  \mbox{ one has }
  \, \, \frac{1}{x_{j}} > \frac{1}{4}\left\|\left(\frac{1}{x_{i}}
  \right)_{i=j}^{n} \right\|_{\frac{H}{k+1-j}} .
$$

We choose $A_{j}=\{j\}$ for $j=1,\dots,k-1$ and
$A_{k}=\{k,\dots,n\}$. Then for every $j\leq k$
$$
  \left\|\left(\frac{1}{x_{i}}\right)_{i\in A_{j}}\right\|_{H}
  \geq \frac{1}{x_{j}} > \frac{1}{4}\left\|\left(\frac{1}{x_{i}}
  \right)_{i=j}^{n}\right\|_{\frac{H}{k+1-j}},
$$
which proves (\ref{Partition3}).

\smallskip

\noindent
{\it Case 3: \ }
$$
  \frac{1}{x_{1}} > \frac{1}{4}\left\|\left(\frac{1}{x_{i}}
  \right)_{i=1}^{n}\right\|_{\frac{H}{k}} \quad \quad
  \mbox{ and there exists } \, \,   j\leq k \, \,  \mbox{ such that }
  \, \, \frac{1}{x_{j}} \leq \frac{1}{4}\left\|\left(\frac{1}{x_{i}}
  \right)_{i=j}^{n} \right\|_{\frac{H}{k+1-j}} .
$$
Let $m$ be the smallest integer such that $m>1$ and
\begin{equation}\label{partition30}
  \frac{1}{x_{m}}\leq\frac{1}{4}\left\|\left(\frac{1}{x_{i}}
  \right)_{i=m}^{n}\right\|_{\frac{H}{k+1-m}} .
\end{equation}
For $1\leq\ell<m$ we choose $A_{\ell}=\{\ell\}$.
Then
$$
   \left\|\left(\frac{1}{x_{i}}\right)_{i\in A_{\ell}}\right\|_{H}
   =\frac{1}{x_{\ell}}>\frac{1}{4}\left\|\left(\frac{1}{x_{i}}
   \right)_{i=\ell}^{n}\right\|_{\frac{H}{k+1-\ell}}
$$
and therefore
$$
\min_{1\leq j< m}
\left\|\left(\frac{1}{x_{i}}\right)_{i=j}^{n}\right\|_{\frac{H}{k-j+1}}
\leq
4 \min_{1\leq j<m}
\left\|\left(\frac{1}{x_{i}}\right)_{i\in A_{j}}\right\|_{H} .
$$
Now we consider the sequence $0<x_{m}\leq x_{m+1}\leq\cdots\leq x_{n}$ and
proceed as in {\it Case 1}. The assumption of {\it Case 1} is fulfilled by 
(\ref{partition30}). The procedure of {\it Case 1}
gives a partition
$A_{m},\dots,A_{k}$ of $\{m,\dots,n\}$ satisfying (\ref{partition31})
$$
4\min_{m\leq \ell\leq k}
\left\|\left(\frac{1}{x_{i}}\right)_{i\in A_{\ell}}\right\|_{H}
\geq
\left\|\left(\frac{1}{x_{i}}\right)_{i=m}^{n}\right\|_{\frac{H}{k+1-m}}.
$$
This completes the proof. 
\kkk
\vskip 2mm

\begin{lem}\label{KminBelow}
Let $p>0$, $1\leq k\leq n$, and $0<x_1\leq x_2\leq ... \leq x_{n}$.
Let $\xi_{1},\dots,\xi_{n}$ be i.i.d. random variables,
$F(t)=\mathbb P \left(\left\{|\xi_{1}|>t\r\}\right)$,
 $N(t)=\ln\frac{1}{F(t)}$,  and
$N_{j}=\frac{2e}{k-j+1}N$, $j=1,\dots,k $.
Then
$$
  \left(1- \tfrac{1}{\sqrt{2 \pi}}\right) \  \max_{1 \leq j \leq k}\
  \| \left(1/x_i\r)_{i=j}^n \| ^{-p}_{N_{j}} \leq
  \mathbb E\operatorname{k-min}_{1\leq i\leq n} |x_{i}\xi_{i}|^p .
$$
\end{lem}

\medskip

\noindent
{\bf Proof.}
Let $c=\left(1- \tfrac{1}{\sqrt{2 \pi}}\right)^{1/p}$.
It is enough to show that for every $k\leq n$ 
\begin{equation}\label{lowgen}
c  \, \| \left(1/x_i\r) _{i=1}^n \| ^{-1}_{N_{1}}
\leq\left(\mathbb E\, \km_{1\leq i\leq n}|x_{i}\xi_{i}|^p\right)^{1/p} .
\end{equation}
Indeed, assume that (\ref{lowgen}) is true. Fix
$j \leq k$. Since
$$
\mathbb E\, \, \km_{1\leq i\leq n}|x_{i} \xi_{i}|^p
\geq \mathbb E\, \, \kjm_{j\leq i\leq n}| x_{i} \xi_{i}|^p ,
$$
(\ref{lowgen}) implies
$$
\left( \mathbb E\, \, \km_{1\leq i\leq n}|x_{i} \xi_{i}|^p\r)^{1/p}
\geq c \, \| \left(1/x_i\r)_{i=j}^n \| ^{-1}_{N_{j}},
$$
for all $1 \leq j \leq k$.

Now we show estimate (\ref{lowgen}).
Fix $\eps >0$ small enough and put
$$
  A=\|\left(1/x_i\r) _{i=1}^n\|_{N_{1}}^{-1} - \eps.
$$
We use  that $1-t\leq -\ln t$ for $t>0$ and that 
$N_{1}=\frac{2e}{k}N=\frac{2e}{k}\ln\frac{1}{F}$ and
we obtain
\begin{eqnarray*}
a&:=& \tfrac{e}{k}\ \sum _{i=1}^{n} G(A /x_i) =
  \tfrac{e}{k}\ \sum _{i=1}^{n}(1-F(A /x_i))   \\
  &\leq&\tfrac{e}{k}\ \sum _{i=1}^{n}\ln\frac{1}{F(A /x_i)}
  =\frac{1}{2}\sum _{i=1}^{n}N_{1}(A /x_{i}) \leq 1/2.
\end{eqnarray*}
Applying Lemma~\ref{forthtwo}, we get
$$
  \mathbb P\left\{ \km_{1\leq i\leq n} |x_{i} \xi_{i}|^p \geq A^p \right\}
  \geq 1 -  \frac{1}{\sqrt{2 \pi k}} \, \, \frac{a ^k}{1-a}
  \geq 1- \frac{1}{ \sqrt{2 \pi}},
$$
as $a\leq 1/2$.
This implies
$$
   \mathbb E\, \, \km_{1\leq i\leq n}|x_{i} \xi_{i}|^p \geq
   A^{p} \ \mathbb P\left\{\km_{1\leq i\leq n}|x_{i} \xi_{i}| ^p
   \geq A^{p} \right\} \geq \left(1- \tfrac{1}{\sqrt{2 \pi}}\right) \ A^p .
$$
Sending $\eps$ to $0$ we obtain the desired result. 
\kkk
\vskip 3mm

\begin{lem}\label{UpperEst1}
Let $p>0$ and $0<x_1\leq x_2\leq ... \leq x_n$
be real numbers. Let $\xi _1, \dots, \xi_n$ be i.i.d. random variables.
Let $F(t)=\mathbb P(|\xi_{1}| > t)$ be strictly decreasing and 
$N=-\ln F$ be an Orlicz function.
Then
$$
\left(1- \tfrac{1}{\sqrt{2 \pi}}\right)
\left\|\left(\frac{1}{x_{i}}\right)_{i=1}^{n}
\right\|_{2eN}^{-p}
\leq\mathbb E \, \, \min_{1\leq i\leq n} |x_{i}\xi _{i}|^{p}
\leq \left(1+\Gamma(1+p)\right)  \ 
 \left\|\left(\frac{1}{x_{i}}\right)_{i=1}^{n}\right\|_{N}^{-p} .
$$
\end{lem}

\medskip

\noindent
{\bf Remark 1. } If $N$ is an Orlicz function then by (\ref{orlcomp}) 
$$
  (2e)^{-p} \|\cdot \| _{N}^{-p} \leq \|\cdot \| _{2eN}^{-p}  . 
$$ 

\noindent
{\bf Remark 2. } The left hand side inequality does not require the condition 
``$N$ is an Orlicz function."

\medskip

\noindent
{\bf Proof.} The left hand inequality follows from Lemma \ref{KminBelow}.

To prove the right hand side inequality we choose
$$
  t_{0}=\left\| \left(\frac{1}{x_{i}}\r) _{i=1}^n\right\|_{N}^{-p} .
$$
Then for  all $t\geq t_{0}$
$$
  \sum_{i=1}^{n}\ln \left( 1/ F\left( t^{1/p}/ x_i\right) \r) \geq 1 .
$$
By (\ref{ExpFunc}) for all
$t\geq t_0$   and all $x_{i}$
$$
(t/t_{0})^{\frac{1}{p}}
\ln\frac{1}{F(t_{0}^{\frac{1}{p}}/x_{i})}
\leq\ln\frac{1}{F(t^{\frac{1}{p}}/x_{i})}.
$$
By Lemma~\ref{lemmin},
$$
 \mathbb P \left\{\min_{1\leq i\leq n}|x_{i} \xi _{i}|^{p}> t
\right\} = \prod_{i=1}^{n} F\left(t^{\frac{1}{p}}/x_{i}\right) =
\exp\left(- \sum_{i=1}^{n}\ln\frac{1}{F\left( t^{\frac{1}{p}}/x_{i}\right)}\right),
$$
and thus for all $t\geq t_{0}$
$$
\mathbb P \left\{\min_{1\leq i\leq n}|x_{i} \xi _{i}|^{p}> t
\right\}
\leq\exp\left(-(t/t_{0})^{\frac{1}{p}}
\sum_{i=1}^{n}\ln\frac{1}{F\left( t_{0}^{\frac{1}{p}}/x_{i}\right)}\right)
\leq \exp\left(-\left(\frac{t}{t_{0}}\right)^\frac{1}{p}\right).
$$
Therefore
\begin{eqnarray*}
\mathbb E \, \, \min_{1\leq i\leq n} |x_{i}\xi _{i}|^{p}
&&=\int_{0}^{\infty}\mathbb P
\left\{\min_{1\leq i\leq n}|x_{i} \xi _{i}|> t^{^{\frac{1}{p}}}
\right\} dt
\\
&&=\int_{0}^{t_{0}}\mathbb P
\left\{\min_{1\leq i\leq n}|x_{i} \xi _{i}|> t^{^{\frac{1}{p}}}
\right\} dt +\int_{t_{0}}^{\infty}\mathbb P
\left\{\min_{1\leq i\leq n}|x_{i} \xi _{i}|> t^{^{\frac{1}{p}}}
\right\} dt
\\
&&\leq t_{0}+\int_{t_{0}}^{\infty}
\exp\left(-\left(\frac{t}{t_{0}}\right)^\frac{1}{p}\right)dt .
\end{eqnarray*}
We substitute $t=t_0 s^p$, then  
$$
\mathbb E \, \, \min_{1\leq i\leq n} |x_{i}\xi _{i}|^{p}
\leq t_{0}+t_{0}p\int_{1}^{\infty}s^{p-1}e^{-s}ds
\leq t_{0}\left(1+p\Gamma(p)\right) ,
$$
which completes the proof. 
\kkk
\vskip 3mm

\begin{lem}\label{KminAbove1}
Let  $1\leq k\leq n$ and $0<x_1\leq x_2\leq ... \leq x_n$.
Let $\xi _1, \dots, \xi_n$ be i.i.d. random variables.
Let $F(t)=\mathbb P(|\xi_{1}| > t)$ be strictly decreasing, 
and let $N=-\ln F$ be an Orlicz function.
Let $M_j = N/(k-j+1)$, $j=1,\dots,k$.
Then
\begin{eqnarray*}
\mathbb E \, \, \km_{1\leq i\leq n} |x_{i}\xi _{i}|
 \leq  8 e\ln (k+1) \ C_N \
  \ \max_{1 \leq j \leq k}\   \| \left(1/x_i\r)
  _{i=j}^n \| ^{-1}_{M_{j}}
 \end{eqnarray*}
where $C_N = \max\{N(1),1/N(1)\}$. 
\end{lem}

\medskip

\noindent{\bf Proof.} The case $k=1$ follows form Lemma~\ref{UpperEst1}. 
We assume $k\geq 2$.

Let $A_1, \ldots A_k$ be the partition of $\{1\ldots n\}$
given by  Lemma~\ref{Partition}. Then for all $q\geq1$
\begin{eqnarray*}
\mathbb E\ \km_{1\leq i\leq n}|x_{i}\xi_{i}|
&\leq&\mathbb E\ \max_{1\leq j\leq k} \min_{i\in A_{j}}|x_{i}\xi_{i}|
\leq\mathbb E\left( \sum_{j=1}^{k} \left|\min_{i\in A_{j}}|x_{i}\xi_{i}|\right|^{q}
\right)^{\frac{1}{q}}   \\
&\leq&\left( \mathbb E\sum_{j=1}^{k} \left|\min_{i\in A_{j}}|x_{i}\xi_{i}|\right|^{q}
\right)^{\frac{1}{q}}
=\left( \sum_{j=1}^{k} \mathbb E\left|\min_{i\in A_{j}}|x_{i}\xi_{i}|\right|^{q}
\right)^{\frac{1}{q}}.
\end{eqnarray*}
By Lemma \ref{UpperEst1} the latter expression is less than 
$$ 
  \left(1+\Gamma(1+q)\right)^{\frac{1}{q}} \  
  \left(\sum_{j=1}^{k}\left\|\left(\frac{1}{x_{i}}\right)_{i\in A_{j}}
  \right\|_{N}^{-q}\right)^{\frac{1}{q}} \leq 2 q k^{1/q} 
  \max_{1\leq j\leq k}\left\|\left(\frac{1}{x_{i}}\right)_{i\in A_{j}}
  \right\|_{N}^{-1}. 
$$
The choice $q=\ln (k+1)$ gives 
$$
  \mathbb E\ \km_{1\leq i\leq n}|x_{i}\xi_{i}| \leq 2 e \ln (k+1) 
  \max_{1\leq j\leq k}\left\|\left(\frac{1}{x_{i}}\right)_{i\in A_{j}}
  \right\|_{N}^{-1}.
$$
By Lemma \ref{Partition}
this expression is smaller than
$$
  8 e \ C_N\ 
\max_{1\leq j\leq k}\left\|\left(\frac{1}{x_{i}}\right)_{i=j}^{n}
\right\|_{M_{j}}^{-1} .
$$
\kkk

\medskip

\noindent
{\bf Proof of Theorem \ref{Kmin1}.}
The lower estimate follows from Lemma \ref{KminBelow}. Since, by 
(\ref{orlcomp}), 
$$
  \|\cdot \| _{N_j} \leq 2 e \|\cdot \| _{M_j},  
$$ 
 the upper estimate follows by Lemma \ref{KminAbove1}. 
\kkk

\section{k-max}

In this section we prove Theorem \ref{TheoK-max}.
We require a result from \cite{GLSW}.
Let $f$ be a random variable with  continuous distribution and
such that $\mbox{\bf E} |f| < \infty$.
Let $t_{n}=t_n (f) =0$, $t_0 = t_0 (f) = \infty$,  and
for $j=1,\dots, n-1$
\begin{equation} \label{zero}
t_{j} = t_{j} (f) = \sup \left\{ t \  | \
\Bbb P\{\omega|\  |f(\omega)| > t \} \geq \tfrac{j}{n} \right\}.
\end{equation}
Since $f$ has the continuous distribution, we have for every $j\geq 1$
$$
 \Bbb P\{\omega|\  |f(\omega)|\geq t_{j}\}=\tfrac{j}{n} .
$$
For  $j=1,\dots,n$ define the sets
\begin{equation} \label{sigma}
\Omega_{j} = \Omega_{j} (f) =\{\omega|\  t_{j}\leq|f(\omega)|< t_{j-1}\} .
\end{equation}
   For all $j=1,\dots,n$ we have
$$
\Omega_{j}=\{\omega|\  t_{j}\leq|f(\omega)|< t_{j-1}\}
=\{\omega|\  t_{j}\leq|f(\omega)|\}
\setminus
\{\omega|\  t_{j-1}\leq|f(\omega)|\} .
$$
Therefore
$$
  \Bbb P(\Omega_{j})=\tfrac{j}{n}-\tfrac{j-1}{n}=\tfrac{1}{n} .
$$
For $j=1,\dots,n$ let
\begin{equation} \label{one}
y_{j} = y_{j} (f) :=  \int_{\Omega_{j}} | f(\omega) | d\Bbb P(\omega) .
\end{equation}
Then
$$
\sum_{j=1}^{n}y_{j}= \mbox{\bf E} |f|
\hskip 10mm
\text{and}
\hskip 10mm
t_{j}\leq ny_{j}<t_{j-1}
\hskip 10mm
\text{for all } \
j=1,\dots,n .
$$
\vskip 3mm

In \cite{GLSW}, Corollary 2  we proved the following
statement.

\begin{lem}\label{cor1}
Let $f_{1},\dots,f_{n}$ be  i.i.d.
random variables such that $\int |f_{i}(\omega)|d\Bbb P(\omega)=1$. Let $M$
be an Orlicz function such that for all $k=1,\dots,n$
$$
 M^{*}\left(\sum_{j=1}^{k}y_{j}\right)=\tfrac{k}{n} .
$$
Then, for all $x\in\mathbb R^{n}$
$$
  c_{1}\|x\|_{M} \leq\int_{\Omega}\max_{1\leq i\leq n}
  |x_{i}f_{i}(\omega)|d\Bbb P(\omega)\leq c_{2}\|x\|_{M},
$$
where $c_1$ and $c_2$ are absolute positive constants.
\end{lem}
\vskip 3mm

This can be reformulated in the following way.

\begin{lem}\label{1max}
Let $\xi_{1},\dots,\xi_{n}$ be  i.i.d.
random variables such that $\int |\xi_{i}(\omega)|d\Bbb P(\omega)=1$. Let $M$
be the Orlicz function such that for all $s\geq0$
$$
  M(s) =\int_{0}^{s}\int_{\frac{1}{t}\leq|\xi_1|}|\xi_1|d\mathbb P dt
  = \int_{\frac{1}{s}\leq|\xi|} \left(s |\xi| - 1\right)  d\mathbb P.
$$
Then, for all $x\in\mathbb R^{n}$ 
$$
  c_{1}\|x\|_{M} \leq\int_{\Omega}\max_{1\leq i\leq n}|x_{i}\xi_{i}(\omega)|
  d\Bbb P(\omega) \leq c_{2}\|x\|_{M},
$$
where $c_1$ and $c_2$ are absolute  positive constants.
\end{lem}
\vskip 3mm

\noindent
{\bf Proof.} By  definition 
$$
   \sum _{i=1}^k y_i = \int _{t_k\leq \xi _1} |\xi _1(\omega )| \ d\mathbb P(\omega)
$$
and
$$
   \mathbb P\left(\left\{ \omega\, \, \mid \, \,
   t_k\leq \xi _1(\omega ) \right\}\right) =\frac{k}{n}.
$$
Therefore the Orlicz function $M^*$ defined by
$$
 M^{*}\left(\int_{t\leq |\xi_{1}|} | \xi_{1}(\omega) | d\mathbb P(\omega)\right)
=\mathbb P\{\omega|t\leq |\xi_{1}(\omega)|\}
$$
satisfies the condition of Lemma~\ref{cor1}.
It is left to observe that the dual function to $M^{*}$ is
$$
M(s) =\int_{0}^{s}\int_{\frac{1}{t}\leq|\xi_1|}|\xi_1|d\mathbb P dt.
$$
This has been verified in Section~\ref{secorl} (see formulae
(\ref{Orlicz1}) and (\ref{OrliczMin})). 
\kkk 
\vskip 3mm

For the next lemma we need the following simple claim.

\begin{cl}\label{simcl}
Let $(x_i)_{i=1}^n$ be a sequence. Then for every $j\leq n-k$
one has
$$
  \operatornamewithlimits{k-max}_{1\leq i\leq n}|x_{i}| \leq
  \operatornamewithlimits{j-min}_{1\leq i\leq k+j-1}|x_{i}|
  + \max_{k+j\leq i\leq n}|x_{i}|.
$$
\end{cl}
\vskip 3mm

\noindent
{\bf Proof.}
If the numbers
$|x_{1}|,\dots, |x_{k+j-1}|$ contain the $k$ biggest of the numbers
$|x_{1}|,\dots,|x_{n}|$, then
$$
 \operatornamewithlimits{j-min}_{1\leq i\leq k+j-1}|x_{i}|
 =\operatornamewithlimits{k-max}_{1\leq i\leq k+j-1}|x_{i}|
 =\operatornamewithlimits{k-max}_{1\leq i\leq n}|x_{i}|.
$$
On the other hand, if the numbers
$|x_{1}|,\dots, |x_{k+j-1}|$ do not contain the $k$ biggest of the numbers
$|x_{1}|,\dots,|x_{n}|$, then at least one of those is contained in
the numbers $|x_{k+j}|,\dots,|x_{n}|$ and therefore
$$
 \max_{k+j\leq i\leq n}|x_{i}| \geq
 \operatornamewithlimits{k-max}_{1\leq i\leq n}|x_{i}| .
$$
\kkk

\begin{lem}\label{ProbEstKmaxAbove}
Let $x_{1}\geq x_{2}\geq\cdots\geq x_{n}>0$.
Let $\xi_{1},\dots,\xi_{n}$ be i.i.d.
random variables and $F(t)=\mathbb P \left(\left\{|\xi_{1}|>t\r\}\right)$.
Suppose that $F$ is strictly decreasing and $N=-\ln F$
is an Orlicz function.
Assume that $M$ is the Orlicz function such that for all $s \geq 0$
$$
 M(s) =\int_{0}^{s}\int_{\frac{1}{t}\leq|\xi|}|\xi|d\mathbb P dt .
$$
Then we have
$$
  \mathbb E
  \operatornamewithlimits{k-max}_{1\leq i\leq n}|x_{i}\xi_{i}|
  \leq c \min_{1\leq j\leq n-k}\Big\{ C_{N}\ln(k+1) \max_{0 \leq \ell \leq j-1}\
  \|\left(1/x_i\r)_{i=1}^{k+\ell}\|^{-1}_{\frac{2e}{\ell+1}N}
  + \|(x_{k+j},\dots,x_{n})\|_{M}\Big\},
$$
where $c_{}$ is an absolute constant
and $C_N= \max\{N(1),1/N(1)\}$.
\end{lem}
\vskip 3mm

\noindent
{\bf Proof.}
 Claim~\ref{simcl} implies
$$
 \mathbb E
 \operatornamewithlimits{k-max}_{1\leq i\leq n}|x_{i}\xi_{i}|   \leq
 \min_{1\leq j\leq n-k}\left(
 \mathbb E\operatornamewithlimits{j-min}_{1\leq i\leq k+j-1}|x_{i}\xi_{i}|
 +\mathbb E\max_{k+j\leq i\leq n}|x_{i}\xi_{i}|\right).
$$
Applying Theorem \ref{Kmin1} to the sequence
$x_{k+j-1} \leq x_{k+j-2} \leq \cdots \leq x_{1}$, we observe that
$$
  \mathbb E\operatornamewithlimits{j-min}_{1\leq i\leq k+j-1}|x_{i}\xi_{i}|
  \leq 16 e^2 C_{N}\ln(k+1)  \  \max_{1 \leq \ell \leq j}\
  \|\left(1/x_i\r)_{i=1}^{k+j-\ell}\|^{-1}_{\frac{2e}{j-\ell+1}N}.
$$
This is the same as
$$
\mathbb E\operatornamewithlimits{j-min}_{1\leq i\leq k+j-1}|x_{i}\xi_{i}|
  \leq 16 e^2 C_{N}\ln(k+1)  \  \max_{0 \leq \ell \leq j-1}\
  \|\left(1/x_i\r)_{i=1}^{k+\ell}\|^{-1}_{\frac{2e}{\ell+1}N}.
$$
On the other hand, Lemma~\ref{1max} implies
$$
\mathbb E\max_{k+j\leq i\leq n}|x_{i}\xi_{i}|
\leq c \ \|(x_{k+j},\dots,x_{n})\|_{M},
$$
where $c$ is an absolute constant. This completes the proof. 
\kkk
\vskip 3mm

\begin{lem}\label{ProbEstKmaxBelow1} \label{pyat}
Let $x_{1}\geq x_{2}\geq\cdots\geq x_{n}>0$. Let
$\xi_{1},\dots,\xi_{n}$ be i.i.d.
random variables and $F(t)=\mathbb{P}(|\xi_{1}| > t)$.
For $k>1$ let
$$
   N_{F,k}(t) = \frac{F\left(1/t \right)}{4(k-1)}  .
$$
Then
$$
\mathbb E
\operatornamewithlimits{k-max}_{1\leq i\leq n}|x_{i}\xi_{i}|
\geq\max\left\{\frac{1}{2}\ \|(x_{k},\dots,x_{n})\|_{N_{F,k}},\ 
\max_{1\leq\ell\leq n-k}\mathbb E
\operatornamewithlimits{\ell-min}_{1\leq i\leq k+\ell-1}
|x_{i}\xi_{i}|\right\} .
$$
In particular,
if  $N(t)=\ln\frac{1}{F(t)}$,
then
$$
  \mathbb E
  \operatornamewithlimits{k-max}_{1\leq i\leq n}|x_{i}\xi_{i}|
  \geq\max\left\{\frac{1}{2}\ \|(x_{k},\dots,x_{n})\|_{N_{F,k}}, \ 
  (1- \tfrac{1}{\sqrt{2 \pi}}) \max_{1\leq\ell\leq n-k}
  \ \max_{1 \leq j \leq \ell}\ \| \left(1/x_i\r) _{i=1}^{k+\ell-j}
  \|_{\frac{2e}{\ell-j+1}N}^{-1} \right\} .
$$
\end{lem}
\vskip 3mm

\noindent
{\bf Proof.}
First we show
$$
  \mathbb E \operatornamewithlimits{k-max}_{1\leq i\leq n}|x_{i}\xi_{i}|
  \geq \frac{1}{2}\ \|(x_{k},\dots,x_{n})\|_{N_{F,k}}.
$$
 We have
$$
  \mathbb P\left\{\omega\left|\operatornamewithlimits{k-max}_{1\leq i\leq n}
  |x_{i}\xi_{i}(\omega)|\leq t\right.\right\} \leq\sum_{j=0}^{k-1} \sum_{A
  \subset \{1, \dots n\} \atop |A|=j} \prod_{i\in A}F\left(\frac{t}{x_{i}}
  \right)\prod_{i\not\in A}\left(1-F\left(\frac{t}{x_{i}}\right)\right) .
$$
Since $x_{1}\geq x_{2}\geq\cdots\geq x_{n}>0$ and $|A^c|\geq n-k+1$,
we observe
$$
  \mathbb P\left\{\omega\left|\operatornamewithlimits{k-max}_{1\leq i\leq n}
  |x_{i}\xi_{i}(\omega)|\leq t\right.\right\}
  \leq\sum_{j=0}^{k-1}\sum_{|A|=j}
  \prod_{i\in A}F\left(\frac{t}{x_{i}}\right)
  \prod_{i=k}^{n}\left(1-F\left(\frac{t}{x_{i}}\right)\right) .
$$
Now we apply the Hardy-Littlewood-Polya inequality (\cite{HLP}), which states that
for non-negative numbers $a_1, \dots, a_m$ one has
$$
  \sum_{A\subset \{1, \dots m\} \atop |A|=j} \prod_{i\in A} a_i \leq
  {m \choose j} \left(\frac{1}{m} \ \sum_{i=1}^{m} a_i
  \right)^{j} \leq \frac{1}{j!} \left(\sum_{i=1}^{m} a_i \right)^{j}.
$$
This implies
$$
  \mathbb P\left\{\omega\left|\operatornamewithlimits{k-max}_{1\leq i\leq n}
  |x_{i}\xi_{i}(\omega)|\leq t\right.\right\} \leq \sum_{j=0}^{k-1}\frac{1}{j!}
  \left(\sum_{i=1}^{n}F\left(\frac{t}{x_{i}}\right)\right)^{j}
  \prod_{i=k}^{n}\left(1-F\left(\frac{t}{x_{i}}\right)\right) .
$$
Since $F\left(\frac{t}{x_{i}}\right)\leq 1$ and
$1-x\leq e^{-x}$ for $x\geq0$, one has
$$
  \mathbb P\left\{\omega\left|\operatornamewithlimits{k-max}_{1\leq i\leq n}
  |x_{i}\xi_{i}(\omega)|\leq t\right.\right\}
  \leq\sum_{j=0}^{k-1}\frac{1}{j!}\left(k-1+\sum_{i=k}^{n}
  F\left(\frac{t}{x_{i}}\right)\right)^{j}
  \exp\left(-\sum_{i=k}^{n}F\left(\frac{t}{x_{i}}\right)\right)  .
$$
Let
$$
 \alpha = \alpha (t) = \frac{1}{k-1} \ \sum_{i=k}^{n}F\left(\frac{t}{x_{i}}\right).
$$
Then
\begin{eqnarray*}
\mathbb P\left\{\omega\left|\operatornamewithlimits{k-max}_{1\leq i\leq n}
|x_{i}\xi_{i}(\omega)|\leq t\right.\right\}
&&\leq e^{-\alpha (k-1)}
\sum_{j=0}^{k-1}\frac{1}{j!}\left((\alpha+1)(k-1)\right)^{j}   \\
&&\leq e^{-\alpha (k-1)}(1+\alpha)^{k-1}
\sum_{j=0}^{k-1}\frac{1}{j!}\left(k-1\right)^{j}   \\
&&\leq e^{-\alpha (k-1)}(1+\alpha)^{k-1}e^{k-1}\\ 
&&= \exp\left(\left(k-1\r)\left(-\alpha +1 + \ln\left(1+\alpha \r) \r) \r).
\end{eqnarray*}
Now put
$$
  t_0 :=  \|(x_{k},\dots,x_{n})\|_{N_{F,k}} \geq 0 . 
$$
If $t_0=0$ we are done. If $t_{0}>0$
then for every $0<\eps<t_0$
$$
 \alpha (t_0 - \eps) = \frac{1}{k-1} \ \sum_{i=k}^{n}F\left(
 \frac{t_0 - \eps }{x_{i}}\right) >4 .
$$
Since $k>1$,  this implies
$$
 \mathbb P\left\{\omega\left|\operatornamewithlimits{k-max}_{1\leq i\leq n}
 |x_{i}\xi_{i}(\omega)|\leq t_0 -\eps \right.\right\}
 \leq\exp((k-1)(-3+\ln 5))\leq 1/2. 
$$
Thus
$$
  \mathbb E
  \operatornamewithlimits{k-max}_{1\leq i\leq n}|x_{i}\xi_{i}|
  \geq  \left(t_0 -\eps\r) \ \mathbb P
  \left\{\omega\left|\operatornamewithlimits{k-max}_{1\leq i\leq n}
  |x_{i}\xi_{i}(\omega)|\geq t_0 -\eps \right.\right\} \geq
  \frac{t_0 -\eps}{2}.
$$
Letting $\eps$ tend to $0$ we obtain the first part of the desired estimate.
\par
Now we show the second part of the estimate.
We observe that for all $l\leq n-k+1$
$$
  \operatornamewithlimits{\ell-min}_{1\leq i\leq k+\ell-1}|x_{i}\xi_{i}(\omega)|
  =\operatornamewithlimits{k-max}_{1\leq i\leq k+\ell-1}|x_{i}\xi_{i}(\omega)|.
$$
This implies
$$
  \mathbb E \operatornamewithlimits{k-max}_{1\leq i\leq n}|x_{i}\xi_{i}|
  \geq \max_{1\leq\ell\leq n-k+1}\mathbb E
  \operatornamewithlimits{\ell-min}_{1\leq i\leq k+\ell-1}|x_{i}\xi_{i}|.
$$

Finally,  the ``In particular" part of Lemma~\ref{pyat} follows from
Lemma~\ref{KminBelow}. Note that in Lemma~\ref{KminBelow} the sequence
$x$ is in increasing order while in Lemma~\ref{pyat} it is in
decreasing order. 
\kkk
\vskip 3mm

In the next lemma we provide a lower estimate on $\|\cdot\| _{N_{F,k}}$,
appearing in Lemma~\ref{pyat}.

\begin{lem}\label{OrlHom}
Let $1< k \leq n$. Let $\xi_{1},\dots,\xi_{n}$ be symmetric i.i.d.
random variables with $\int |\xi_{i}(\omega)|d\mathbb P(\omega)=1$.
Let $F(t)=\mathbb P(|\xi_{1}|>t)$ be a strictly decreasing function 
such that $N=-\ln F$ is an Orlicz function. 
Let $M$ be the Orlicz function defined by 
$$
  M(s)=\int_{0}^{s}\int_{\frac{1}{t}\leq|\xi_1|}|\xi_1|d\mathbb P dt.
$$
Let
$$
N_{F,k}(t) = \frac{F\left(1/t \right)}{4(k-1)} .
$$
Let $k_{0}=\left[\frac{4(k-1)}{F(1)}\right]$ and assume that $k+k_0\leq n$.
Then for all $x_{1}\geq\cdots\geq x_{n}>0$ 
$$
  \|(x_{k+k_{0}},\dots,x_{n})\|_{M}
\leq \left( 1 + \frac{\ln (8(k-1))}{N(1)} \r)\ 
  \|(x_{k},\dots,x_{n})\|_{N_{F,k}} . 
$$

\end{lem}
\vskip 3mm

\noindent
{\bf Remark. } The condition {\it $N=-\ln F$ is an Orlicz function} in this Lemma 
can be substituted with the condition {\it there is a constant $c_{2}\geq 1$ such that 
for all $s \in (0,1/c_2]$ and $t\in (0, 1/(4 c^2_{2}) ]$
\begin{equation}\label{SubMultEst2}
F^{-1}(s)F^{-1}(t)\geq F^{-1}(c_{2}st)
\end{equation}
and such that for all
$t\geq F^{-1}(\frac{1}{c_{2}})$
we have
\begin{equation}\label{semna}
\int_{t\leq|\xi_{1}|}|\xi_{1}|d\mathbb P\leq c_{2} t F(t) .
\end{equation}
}
Then the conclusion will be 
$$
  \|(x_{k+k_{0}},\dots,x_{n})\|_{M}\leq F^{-1}\left(
  \alpha \right)\
  \|(x_{k},\dots,x_{n})\|_{N_{F,k}} , 
$$
where $\alpha = 1/(4 c_2^2 (k-1))$.

\vskip 3mm

\noindent
{\bf Proof.}
Since both functions $\|\cdot\|_{M}$ and
$\|\cdot \|_{N_{F,k}}$ are homogeneous, we may assume
that $\|(x_{k+k_{0}},\dots,x_{n})\|_{M}=1$. Thus, without loss of
generality, we can assume that
$\sum_{i=k+k_{0}}^{n}M(x_{i})=1$
(otherwise we pass from the sequence $\{x_i\}_i$ to $\{x_i/(1+\eps )\}_i$
for an suitably small $\eps>0$).
\par
We put 
$$
  A:= F^{-1} (\alpha) = 1 + \frac{\ln (8(k-1))}{N(1)} .
$$
Note that by (\ref{conin}), $N(A) \geq A N(1)\geq \ln 8 >2$.  

\smallskip

\noindent
{\bf Case 1:} $x_{k+k_0}\geq 1/ A$.
Then $x_k\geq x_{k+1}\geq \ldots \geq x_{k+k_0} \geq 1/A$.

\medskip

Since $F$ is a
decreasing function, $1/F^{-1}$ is increasing and
$$
\sum_{i=k}^{n}
F\left(\left(x_{i} A \right)^{-1}\right)
\geq \sum_{i=k}^{k+k_{0}}
F\left(\left(x_{i} A \right)^{-1}\right)
\geq (k_{0}+1) F(1)> 4(k-1) .
$$
This means that
$$
  \|(x_{k},\dots,x_{n})\|_{N_{F,k}} \geq 1/A . 
$$

\medskip

\noindent
{\bf Case 2:}
$x_{k+k_{0}} < 1/A$. 
Then $1/A > x_{k+k_0}\geq \ldots \geq x_{n}$.

\smallskip

Since
$\int_{\frac{1}{t}\leq |\xi_{1}|}|\xi_{1}|d\mathbb P$ is an
increasing function of $t$, we observe
$$
  M(s)=\int_{0}^{s}\int_{\frac{1}{t}\leq |\xi_{1}|}|\xi_{1}|d\mathbb P dt
  \leq s\ \int_{\frac{1}{s}\leq |\xi_{1}|}|\xi_{1}|d\mathbb P .
$$
By (\ref{SubMultEst3}), applied with $t=1/s$, we obtain that for all 
positive $s$ 
$$
  \int_{1/s \leq |\xi_{1}|}|\xi_{1}|d\mathbb P\leq 
   \left(1+ \frac{1}{N(1/s)} \r)\ \frac{1}{s} \ F(1/s) .
$$
Recall that $N$ is increasing and $N(A)>2$. Thus for all $s\leq 1/A$ 
we have 
$$
  M(s)\leq s\ \int_{\frac{1}{s}\leq |\xi_{1}|}|\xi_{1}|d\mathbb P
  \leq 2 \  F(1/s).
$$
By the condition of Case 2, $x_i \leq 1/A$ for $i\geq k+k_0$. This implies 
\begin{equation}\label{est30}
  1=\sum_{i=k+k_{0}}^{n}M(x_{i}) \leq 2 \
  \sum_{i=k+k_{0}}^{n}F\left(\frac{1}{x_{i}}\right).
\end{equation} 

Now, by (\ref{conin}), we have $N(y) \geq \beta N(y/\beta)$
for every $y\geq 0$ and $\beta \geq 1$. Since $N=-\ln F$, we observe 
$$
  F(y) \leq F(y/\beta) ^{\beta} 
$$
for every $y\geq 0$ and $\beta \geq 1$. Since $F$ is decreasing, it implies 
$$
 F(y) \leq F(y/\beta)  \ F(1) ^{\beta - 1} 
$$
for every $y\geq \beta \geq 1$. Applying the last inequality 
with $y=1/x_i$  and $\beta = A$, we obtain for every 
 $i\geq k+k_0$
$$
  F(1/x_i) \leq F (1/(A x_i)) \ F(1) ^{A - 1} .
$$ 
By (\ref{est30}), 
$$ 
 \sum_{i=k}^{n} F\left(1/ (A x_i) \right) \geq \frac{1}{F(1) ^{A - 1}} 
 \, \sum_{i=k+k_0}^{n} F\left(1/ x_i \right) \geq \frac{1}{2 F(1) ^{A -1}}.  
$$ 
Now, by the choice of $A$,  
$$
 A -1 = \frac{\ln (8(k-1))}{\ln(1/F(1))}, 
$$
and hence
$$
  2 F(1) ^{A -1} = \frac{1}{4(k-1)}. 
$$
Thus, 
$$
 \sum_{i=k}^{n} F\left(1/ (A x_i) \right) \geq \frac{1}{4(k-1)}, 
$$
which implies 
$$
  \|(x_{k},\dots,x_{n})\|_{N_{F,k}}\geq 1/A.
$$
It completes the proof.
\kkk

\vskip 3mm

Now we  prove Theorem \ref{TheoK-max}.

\medskip

\noindent
{\bf Proof of Theorem \ref{TheoK-max}.}
Let $k_{0}=\left[\frac{4(k-1)}{F(1)}\right]$.
First we show the right hand side inequality. With $j=k_{0}$
in Lemma \ref{ProbEstKmaxAbove}
$$
   \mathbb E
  \operatornamewithlimits{k-max}_{1\leq i\leq n}|x_{i}\xi_{i}|
  \leq c\ \min_{1\leq j\leq n-k}\Big\{ C_{N}\ln(k+1) \max_{0 \leq \ell \leq j-1}\
  \|\left(1/x_i\r)_{i=1}^{k+\ell}\|^{-1}_{\frac{2e}{\ell+1}N}
  + \ \|(x_{k+j},\dots,x_{n})\|_{M}\Big\} .
$$
To show the left hand side inequality we apply
Lemma \ref{ProbEstKmaxBelow1} with $l=k_{0}$:
$$
  \mathbb E
  \operatornamewithlimits{k-max}_{1\leq i\leq n}|x_{i}\xi_{i}|
  \geq\max\left\{\frac{1}{2}\ \|(x_{k},\dots,x_{n})\|_{N_{F,k}} ,
  \left(1- \tfrac{1}{\sqrt{2 \pi}}\right) \ \max_{1 \leq j \leq k_{0}}\
  \| \left(1/x_i\r) _{i=1}^{k+k_{0}-j} \|_{\frac{2e}{k_{0}-j+1}N}^{-1} \right\} .
$$
By Lemma \ref{OrlHom}
we have for all $x$ with $x_{1}\geq\cdots\geq x_{n}>0$
$$
  \|(x_{k+k_{0}},\dots,x_{n})\|_{M} \leq 
 A \ \|(x_{k},\dots,x_{n})\|_{N_{F,k}}  ,
$$
where $A= 1 + \frac{\ln (8(k-1))}{N(1)}$. 
Thus
\begin{eqnarray*}
  \mathbb E \operatornamewithlimits{k-max}_{1\leq i\leq n}|x_{i}\xi_{i}|
& \geq & \frac{1}{4 A} \ 
\|(x_{k+k_{0}},
  \dots,x_{n})\|_{M} + \frac{1}{2}\left(1- \tfrac{1}{\sqrt{2 \pi}}\right)
  \ \max_{1 \leq j \leq k_{0}}\   \| \left(1/x_i\r) _{i=1}^{k+k_{0}-j}
  \|_{\frac{2e}{k_{0}-j+1}N}^{-1}
\\
& = & \frac{1}{4 A}\ 
\|(x_{k+k_{0}},
  \dots,x_{n})\|_{M} + \frac{1}{2}\left(1- \tfrac{1}{\sqrt{2 \pi}}\right)
  \ \max_{0 \leq \ell \leq k_{0}-1}\   \| \left(1/x_i\r) _{i=1}^{k+\ell}
  \|_{\frac{2e}{\ell+1}N}^{-1}.
\end{eqnarray*}
\kkk

\vskip 3mm

\section{The Gaussian case}
\label{gausscl}

In this section we verify that Gaussian $N(0, 1)$-variables
satisfy the hypotheses of Theorems \ref{Kmin1} and \ref{TheoK-max}.
Hence in this section we consider only
standard Gaussian variables and we denote them by $\xi$.
Accordingly,
$$
  F(t)=\mathbb P \left(\left\{| \xi |> t\r\}\right)=
  \sqrt{\tfrac{2}{\pi}} \int_{t}^{\infty}e^{-\frac{s^{2}}{2}}ds
$$
for $t\geq 0$
and $N=\ln\frac{1}{F} = - \ln F$. We check the properties
of the functions $F$ and $N$ in several claims.
\vskip 3mm

\begin{cl} For every $t > 0$  and $A>0$ 
 \begin{equation}\label{gadi}
   F(t) =\sqrt{\tfrac{2}{\pi}} \int_{t}^{\infty}e^{-\frac{s^{2}}{2}}ds
  \leq \sqrt{\tfrac{2}{\pi}}\ \frac{1}{t} \ \exp\left(-\frac{t^{2}}{2}\r)
\end{equation}
and
 \begin{equation}\label{gaditwo}
  \sqrt{\tfrac{2}{\pi}}\  A \ \exp\left(-\left(t+A\r)^2/2\r) \leq F(t).
\end{equation}
In particular,  for $A=1/t$
 \begin{equation}\label{gaditri}
  \sqrt{\tfrac{2}{\pi}}\  \frac{1}{e t} \ \exp\left(-\left(t^{2}+1/t^2\r)/2\r)
  \leq F(t) .
\end{equation}

\end{cl}

\medskip

\noindent{\bf Proof. }
Both estimates follow immediately by integration. Indeed,
the upper estimate follows from
$$
  t \int_{t}^{\infty}e^{-\frac{s^{2}}{2}}ds \leq \int_{t}^{\infty}
  s e^{-\frac{s^{2}}{2}}ds = e^{-\frac{t^{2}}{2}} ,
$$
while the lower estimate follows from
$$
   \int_{t}^{\infty}e^{-\frac{s^{2}}{2}}ds \geq  \int_{t}^{t + A}
   e^{-\frac{s^{2}}{2}}ds \geq A\exp\left(- \left(t + A \r)^2/2 \r).
$$
\kkk

\begin{cl} \label{firstthea}
  $N$ is an Orlicz function.
\end{cl}

\vskip 3mm

\noindent{\bf Proof. }
$N$ is an increasing function on $[0, \infty)$
such that $N(t)=0$ if and only if $t=0$. We have to show that
$N$ is convex. To do so,  we show that $N^{\prime \prime}(t) \geq 0$
for $t\geq 0$.
$$
   N^{\prime \prime}(t)
   =\left(\frac{e^{-\frac{t^{2}}{2}}}
   {\int_{t}^{\infty}e^{-\frac{s^{2}}{2}}ds}\right)^{\prime}
   =\frac{- t e^{-\frac{t^{2}}{2}} \int_{t}^{\infty}
   e^{-\frac{s^{2}}{2}}ds + e^{-t^{2}}  }{ \left(
   \int_{t}^{\infty} e^{-\frac{s^{2}}{2}} ds \r)^{2} } =
  \frac{ e^{-\frac{t^{2}}{2}} \left( e^{-\frac{t^{2}}{2}} - t
   \int_{t}^{\infty} e^{-\frac{s^{2}}{2}}ds \r)}{ \left(
   \int_{t}^{\infty} e^{-\frac{s^{2}}{2}}ds\r)^{2}} ,
$$
which is non-negative by  (\ref{gadi}).
\kkk

Claim \ref{firstthea} shows that $F$ and $N$ satisfy
the hypothesis of Theorem~\ref{Kmin1} and \ref{TheoK-max}. Next we 
estimate the Orlicz norm $N_j$.

\begin{cl} \label{firstthec}
Let $N=-\ln F$,
$$
  H (t) =   \left\{
  \begin{array}{ll}
   t  \quad  \quad \quad  \quad \mbox{ for } \
   0\leq  t < 1
   \\
   t^2  \, \, \,     \quad \quad  \quad
  \mbox{ for } \  t \geq 1.
   \end{array}
\right.
$$
Then $H$ is an Orlicz function and for every $t\geq 0$
$$
     \left(2 \pi e\right)^{-1/2}\  H(t) \leq N(t) \leq 4.5 \  H(t) .
$$
In particular, if $k\leq n$ and $N_j$, $j\leq k$, as in (\ref{OrlNj})
then for every $t\geq 0$
$$
  \sqrt{\frac{2 e}{\pi}}\ \frac{1}{k-j+1} \ H(t) \leq N_j (t)
  \leq \frac{9 e}{k-j+1}\ H(t).
$$
\end{cl}

\medskip

\noindent{\bf Proof. } Clearly, $H$ is an Orlicz function.

For every $0\leq t\leq \sqrt{\pi/8}$ we have
$$
  \frac{1}{2} \leq 1-\sqrt{\tfrac{2}{\pi}}\ t \leq F(t) =
  1 - \sqrt{\tfrac{2}{\pi}} \int_{0}^{t}e^{-\frac{s^{2}}{2}}ds
  \leq 1-\sqrt{\tfrac{2}{e\pi}}\ t.
$$
Since $(x-1)/2 \leq \ln x \leq x-1$ on $[1, 2]$, we
observe for $0\leq t\leq \sqrt{\pi/8}$
$$
  N(t) = \ln\frac{1}{F(t)} \leq\frac{1}{F(t)}-1 \leq \frac{1}{1
  -\sqrt{\tfrac{2}{\pi}} \ t} - 1 \leq \frac{\sqrt{\tfrac{2}{\pi}}
  \ t }{1-\sqrt{\tfrac{2}{\pi}}\ t} \leq \sqrt{\tfrac{8}{\pi}} \ t
$$
 and
$$
  N(t) = \ln\frac{1}{F(t)} \geq \frac{1}{2}\left(\frac{1}{F(t)}
  -1\right) \geq \frac{1}{2}\left(\frac{1}{1-\sqrt{\tfrac{2}{e \pi}}
  \ t}-1\right) \geq \frac{t}{\sqrt{ 2 e \pi} } .
$$
This shows the desired result for $0\leq t\leq \sqrt{\pi/8}$.

Consider now  the function $f(t)= N(t) - t^2/2$ and  observe that $f(0)=0$.
By (\ref{gadi}) we have $f'(t)\geq 0$ for $t\geq 0$. Thus, for every
$t\geq 0$ one has $N(t)\geq t^2/2$.

Finally, applying (\ref{gaditwo}) with $A=\sqrt{\pi/2}$, we have for
$t\geq \sqrt{\pi/8}$ (then $t+A\leq 3t$) that
$$
   F(t) \geq \exp\left(-9 t^2/2\r).
$$
This  implies $t^2/2 \leq N(t) \leq 9 t^2/2$ for $t\geq \sqrt{\pi/8}$.
In particular, $H/\sqrt{2\pi e} \leq N \leq 9 H/2$.
\kkk

\smallskip

\noindent
Y. Gordon, {\small Dept. of Math.}, {\small Technion},
{\small Haifa 32000, Israel}, {\small \tt gordon@techunix.technion.ac.il}

\smallskip

\noindent
A. E. Litvak,  {\small Dept. of Math. and Stat. Sciences},
{\small University of Alberta}, {\small Edmonton, AB, Canada T6G 2G1},
{\small \tt  alexandr@math.ualberta.ca}

\smallskip

\noindent
C. Sch\"utt, {\small Mathematisches Seminar},
{\small Christian Albrechts Universit\"at}, {\small 24098 Kiel, Germany},
\newline
{\small \tt schuett@math.uni-kiel.de}

\smallskip

\noindent
E. Werner, {\small Dept. of Math.},
{\small Case Western Reserve University}, {\small Cleveland, Ohio 44106,
U.S.A.}  and {\small Universit\'{e} de Lille 1}, {\small UFR de Math\'{e}matique},
{\small 59655 Villeneuve d'Ascq, France},
{\small \tt emw2@po.cwru.edu}


\bigskip



\begin{thebibliography}{1234}
{\footnotesize

\bibitem{AN}
{\sc B. C. Arnold, N. Narayanaswamy}, {\em Relations, Bounds and
Approximations for Order Statistics}, Lecture Notes in Statistics,
53, Berlin etc.: Springer-Verlag. viii (1989).

\bibitem{BCh} {\sc N.
Balakrishnan, W.W.S. Chen},
 {\em Handbook of Tables for order
Statistics from lognormal Distributions with Applications},
Amsterdam, Netherlands, Kluwer Academic Publishers (1999).

\bibitem{BCo} {\sc N. Balakrishnan, A.C.
Cohen}, {\em Order Statistics and Inference},  New York, NY:
Academic Press (1991).


\bibitem{BaDaDeVWa}
{\sc R Baraniuk, M Davenport, R DeVore, M Wakin}, {\em The
Johnson-Lindenstrauss lemma meets compressed sensing}, preprint.

\bibitem{BCS}{\sc N. Balakrishnan, E. Castillo and  J.M. Sarabia,
editors}, {\em Advances in Distribution Theory, Order Statistics,
and Inference Series: Statistics for Industry and Technology},
 Birkh\"auser (2006).

\bibitem{BL}
{\sc Y. Benjamini, M. Leshno}, {\em Statistical Methods for Data
Mining}, Data Mining and Knowledge Discovery Handbook, Springer US
(2005).

\bibitem{BR1}{\sc
N. Balakrishnan and C. R. Rao, editors}, {\em Handbook of Statistics
16: Order Statistics: Theory and Methods Elsevier}, Amsterdam
(1999a).

\bibitem{BR2}{\sc N.
Balakrishnan and C. R. Rao, editors}, {\em Handbook of Statistics
16: Order Statistics: Applications},  Elsevier, Amsterdam (1999b).

\bibitem{CRT}
{\sc E. J. Candes, J. Romberg, T. Tao}, {\em Robust Uncertaintity
Principles: Exact Signal Reconstruction from Highly Incomplete
Frequency Information}, IEEE Trans. Inf. Theory, to appear.

\bibitem{CoDDeV}
{\sc A. Cohen, W. Dahmen and R. DeVore}, {\em Compressed Sensing and
Best $k$-term Approximation}, preprint.

\bibitem{DN}
{\sc  H. A. David,  H. N. Nagaraja}, {\em Order statistics}, 3rd
ed., Wiley Series in Probability and Statistics. Chichester: John
Wiley \& Sons (2003).

\bibitem{DG}
{\sc A. Dimitriyuk, Y. Gordon}, {\em Generalizing the
Johnson-Lindenstrauss lemma to $k$-dimensional affine subspaces},
preprint.

\bibitem{D}
{\sc D. Donoho}, {\em Compressed Sensing}, IEEE Trans. Information
Theory, 52 (2006), 1289--1306.

\bibitem{GK} 
{\sc E.D. Gluskin and S. Kwapie\'n}  
{\em Tail and moment estimates for sums of independent random variables}, 
Stud. Math. 114 (1995), 303-309. 

\bibitem{Go2}
{\sc Y. Gordon}, {\em Majorization of Gaussian processes and
Geometric Applications}, Prob. Th. Rel. Fields, 91, No. 2 (1992),
251--267.

\bibitem{GGMP}
{\sc Y. Gordon, O. Gu\'edon, M. Meyer and A. Pajor}, {\em On the
Euclidean sections of some Banach spaces and operator spaces}, Math.
Scandinavica 91 (2002), 247--268.

\bibitem{GLMP} {\sc Y. Gordon, A. E. Litvak,  S. Mendelson, A. Pajor}, {\em Gaussian averages
of interpolated bodies and applications to approximate
reconstruction}, J. Approx. Theory,  149 (2007), 59--73.

\bibitem{GLSW}
{\sc Y. Gordon, A. E. Litvak, C. Sch\"utt, E. Werner}, {\em Orlicz
Norms of Sequences of Random Variables}, Ann. of Prob., 30 (2002),
1833--1853.

\bibitem{GLSW2}
{\sc Y. Gordon, A. E. Litvak, C. Sch\"utt, E. Werner},
{\em Geometry of spaces between zonoids and polytopes},
Bull. Sci. Math.,  126 (2002), 733--762.

\bibitem{GLSW3}
{\sc Y. Gordon, A. E. Litvak, C. Sch\"utt, E. Werner},
{\em On the minimum of several  random variables},
Proc. Am. Math. Soc. 134, No. 12, 3665--3675 (2006). 

\bibitem{GLSW4}
{\sc Y. Gordon, A. E. Litvak, C. Sch\"utt, E. Werner},
{\em Minima of sequences of Gaussian random variables},
C.R. Acad. Sci. Paris, Ser 1, Math.,
340 (2005), 445--448.

\bibitem{Gu}
{\sc O. Gu\`{e}don}, {\em Gaussian Version of a Theorem of Milman
and Schechtman}, Positivity 1, No.1 (1997), 1--5.


\bibitem{HLP}
{\sc G. H. Hardy, J. E. Littlewood and G. Polya}, {\em
Inequalities}, 2nd ed.,  Cambridge, The University Press. XII
(1952).

\bibitem{JL}
{\sc W. B. Johnson, J. Lindenstrauss}, {\em Extensions of Lipschitz
Mappings into a Hilbert Space}, Contemp. Math. 26 (1984), 189--206.


\bibitem{J}
{\sc M. Junge}, {\em The optimal order for the $p$-th moment of sums of 
independent random variables with respect to symmetric norms and related 
combinatorial estimates}, Positivity 10 (2006), 201-230. 

\bibitem{KrRu}
{\sc M.A. Krasnosel'skii and Ya. B. Rutickii}, {\em Convex Functions
and Orlicz Functions}, P. Noordhoff, Groningen (1961).


\bibitem{KS1} {\sc S. Kwapien, C. Sch\"utt},
{\em Some combinatorial and probabilistic inequalities
and their application to Banach space theory},
Studia Math. 82 (1985), 91--106.


\bibitem{KS2} {\sc S. Kwapien, C. Sch\"utt},
{\em Some combinatorial and probabilistic inequalities and their
application to Banach space theory. II},
Studia Math. 95 (1989),  141--154.


\bibitem{Lat} {\sc R. Lata{\l}a}, preprint. 


\bibitem{LT}
{\sc J. Lindenstrauss and L. Tzafriri}, {\em Classical Banach 
Spaces I}, Springer-Verlag (1977).


\bibitem{MZ} {\sc S. Mallat and O. Zeitouni}, {\em Optimality of the 
Karhunen-Loeve basis in nonlinear reconstruction}, preprint. 

\bibitem{MaPa} {\sc V. A. Marchenko and  L. A. Pastur}, {\em Distribution
  of eigenvalues in certain sets of random matrices},
  Mat. Sb. (N.S.), 72 (1967), 407--535 (Russian).

\bibitem{MPT}
{\sc S. Mendelson, A. Pajor, N. Tomczak-Jaegermann}, {\em
Reconstruction and Subgaussian Operators in Asymptotic Geometric
Analysis}, GAFA, Geom. and Funct. Anal., 17 (2007), 1248--1282.

\bibitem{MS}
{\sc S. Montgomery-Smith}, {\em Rearrangement invariant norms of 
symmetric sequence norms of independent sequences of random variables},  
Isr. J. Math. 131 (2002), 51--60. 

\bibitem{RaRe}
{\sc M.M. Rao und Z.D. Ren}, {\em Theory of Orlicz Spaces}, Marcel
Dekker (1991).

\bibitem{RU} {\sc M. Rudelson} 
{\em Lower estimates for the singular values of random matrices},  
C. R. Math. Acad. Sci. Paris  342  (2006),  no. 4, 247--252. 

\bibitem{SZ1} {\sc S. J. Szarek}, {\em Spaces with large distance to 
$\ell ^n_{\infty}$ and random matrices},   
Amer. J. Math.  112  (1990),  no. 6, 899--942. 

\bibitem{SZ2} {\sc S. J. Szarek}, 
{\em Condition numbers of random matrices},   
J. Complexity  7  (1991),  no. 2, 131--149. 

\bibitem{W1}
{\sc L. Wasserman}, {\em All of Statistics: A Concise Course in
Statistical Inference}, Springer Texts in Stat.  (2004).

}

\end{thebibliography}
\end{document}